\documentclass{amsart}

\usepackage{amsmath}
\usepackage{amscd}
\usepackage{amssymb}

\newcommand{\cal}{\mathcal}

\newcommand{\bE}{{\Bbb E}}

\newcommand{\bP}{{\Bbb P}}

\newcommand{\cC}{{\cal C}}

\newcommand{\cM}{{\cal M}}

\newcommand{\Mbar}{\overline{\cM}}

\DeclareMathOperator{\Aut}{Aut} \DeclareMathOperator{\ch}{ch}

\newtheorem{theorem}{Theorem}[section]
\newtheorem{corollary}{Corollary}[section]

\newtheorem{lemma}[theorem]{Lemma}

\theoremstyle{remark}

\theoremstyle{definition}

\begin{document}

\title{Mari\~no-Vafa Formula and  Hodge Integral Identities}
\author{Chiu-Chu Melissa Liu}
\address{Department of Mathematics, Harvard University, 
Cambridge, MA 02138, USA \and 
Center of Mathematical Sciences, Zhejiang University, 
Hangzhou, Zhejiang 310027, China}
\email{ccliu@math.harvard.edu}
\author{Kefeng Liu}
\address{Center of Mathematical Sciences, Zhejiang University,
Hangzhou, Zhejiang 310027, China \and
Department of Mathematics, University of California at Los Angeles,
Los Angeles, CA 90095-1555, USA}
\email{liu@cms.zju.edu.cn, liu@math.ucla.edu}
\author{Jian Zhou}
\address{Department of Mathematical Sciences, Tsinghua University,
Beijing 100084, China \and 
Center of Mathematical Sciences, Zhejiang University,
Hangzhou, Zhejiang 310027, China}

\email{jzhou@math.tsinghua.edu.cn}

\begin{abstract}

We derive some Hodge integral identities
by taking various limits of the Mari\~no-Vafa formula
using the cut-and-join equation.
These identities include the formula of general
$\lambda_g$-integrals, the formula of $\lambda_{g-1}$-integrals
on $\Mbar_{g,1}$, the formula of cubic $\lambda$ integrals on $\Mbar_g$,
and the ELSV formula relating  Hurwitz numbers and Hodge integrals. In particular, our proof
of the MV formula by the cut-and-join equation leads to
a new and simple proof of the $\lambda_g$ conjecture. 
We also present a proof of the ELSV formula completely 
parallel to our proof of the Mari\~{n}o-Vafa formula.
\end{abstract}

\maketitle

\section{Introduction}
Based on duality between Gromov-Witten theory on Calabi-Yau threefolds
and Chern-Simons gauge theory on three manifolds, 
Mari\~no and Vafa \cite{Mar-Vaf} conjectured a closed formula
on certain Hodge integrals  in terms of representations of symmetric
groups. In \cite{LLZ1, LLZ}, the authors of the present paper
proved the Mari\~{n}o-Vafa formula by showing that both sides
of the formula satisfy  the same differential
equation (the cut-and-join equation) and have the same
initial values. For a different approach see \cite{Oko-Pan}.

The first goal of this note is to explain close relationships
between the MV formula and the ELSV
formula relating Hurwitz numbers $H_{g,\mu}$ to Hodge integrals:
\begin{equation}\label{eqn:elsv}
H_{g, \mu} = \frac{(2g-2+|\mu|+l(\mu))!}{|\Aut (\mu)|}
\prod_{i=1}^{l(\mu)} \frac{\mu_i^{\mu_i}}{\mu_i!}
\int_{\Mbar_{g, l(\mu)}} \frac{\Lambda_g^{\vee}(1)}{\prod_{i=1}^{l(\mu)} (1 - \mu_i \psi_i)}.
\end{equation}
The ELSV formula was first derived by Ekedahl, Lando, Shapiro,
Vainshtein \cite{ELSV}.
By the Burnside formula, Hurwitz numbers are related to 
representations of symmetric groups, hence so is the Hodge integral
on the right hand side of (\ref{eqn:elsv}).
In this paper, we will present a proof of the
ELSV formula completely parallel to our proof of the MV
formula \cite{LLZ1, LLZ}, and compare this proof with proofs given in \cite{Gra-Vak}
and \cite[Section 7]{LLZ}. We will also explain how to obtain the ELSV formula
by taking a particular limit of the MV formula, so the  ELSV formula
can also be viewed as an immediate corollary of the MV formula.

The second goal of this note is to show that the following formula
of $\lambda_g$-integrals can be easily derived from  the MV formula
and the cut-and-join equation. 
\begin{equation}\label{eqn:formula}
\int_{\Mbar_{g, n}} \psi_1^{k_1} \cdots \psi_n^{k_n}\lambda_g =
\begin{pmatrix} 2g+n-3 \\ k_1, \dots, k_n\end{pmatrix}
\frac{2^{2g-1}- 1}{2^{2g-1}} \frac{|B_{2g}|}{(2g)!}
\end{equation}
where $k_1 + \cdots + k_n = 2g-3+n$, and $B_{2g}$ are Bernoulli numbers.
Let
\begin{equation}\label{eqn:bg}
b_g=\left\{\begin{array}{ll}
1,& g=0,\\
\int_{\Mbar_{g,1}}\psi_1^{2g-2}\lambda_g, & g>0. 
\end{array}\right.
\end{equation}
Formula (\ref{eqn:formula}) also follows from 
\begin{equation}\label{eqn:sin}
\sum_{g\geq 0} b_g t^{2g}=\frac{t/2}{\sin(t/2)}
=1+\sum_{g\geq 1} \frac{2^{2g-1}-1}{2^{2g-1}}\frac{|B_{2g}|}{(2g)!}t^{2g}
\end{equation}
proved in \cite{Fab-Pan1}, 
and the $\lambda_g$-conjecture 
\begin{equation}\label{eqn:lamg}
\int_{\Mbar_{g, n}} \psi_1^{k_1} \cdots \psi_n^{k_n}\lambda_g =
\begin{pmatrix} 2g+n-3 \\ k_1, \dots, k_n\end{pmatrix}b_g.
\end{equation}
The $\lambda_g$-conjecture (\ref{eqn:lamg}) was found in \cite{Get-P} as a 
consequence of the degree 0 Virasoro conjecture of $\bP^1$ and 
was first proved in \cite{Fab-Pan1}. Later the Virasoro conjecture 
was proved for projective spaces \cite{Giv} and curves 
\cite{Oko-Pan1}; both cases include $\bP^1$ as a special case.

Our proof of the MV formula by the cut-and-join equation 
relies only on the formula (\ref{eqn:sin}) of $b_g$ ($\lambda_g$-integrals on $\Mbar_{g,1}$),   
so we obtain a new and simple proof of the $\lambda_g$ conjecture.
Moreover, the ELSV formula and the formula (\ref{eqn:formula}) of 
$\lambda_g$-integrals are just two particular limits of the MV
formula. On the other hand, another proof
of the MV formula by bilinear localization equations
\cite{Oko-Pan} relies on the ELSV formula and the formula 
(\ref{eqn:formula}) of general $\lambda_g$-integrals.

Finally, we will show that  the following identities of Hodge integrals proved in 
\cite{Fab-Pan1} are also easy consequences of the Mari\~{n}o-Vafa formula and
the cut-and-join equation.
\begin{equation}\label{eqn:cubic}
\int_{\Mbar_{g}} \lambda_{g-2}\lambda_{g-1}\lambda_g =
\frac{1}{2(2g-2)!} \frac{|B_{2g-2}|}{2g-2} \frac{|B_{2g}|}{2g},
\end{equation}
\begin{equation}\label{eqn:MinusOne}
\int_{\Mbar_{g, 1}} \frac{\lambda_{g-1}}{1-\psi_1}
= b_g \sum_{i=1}^{2g-1} \frac{1}{i}
- \frac{1}{2} \sum_{\substack{g_1+g_2=g\\g_1,g_2 > 0}}
\frac{(2g_1-1)!(2g_2-1)!}{(2g-1)!}b_{g_1}b_{g_2}.
\end{equation}

The rest of this paper is organized as follows. In Section  \ref{sec:munford},
we derive some identities from the Mumford's relations which
will be used in later derivations. In Section \ref{sec:MVELSV}, we recall the precise
statements of the MV formula and the ELSV formula, and establish
the relationships between them described above. In Section \ref{sec:MVlambda},
we derive the formula (\ref{eqn:formula}) of general $\lambda_g$-integrals
from the MV formula and the cut-and-join equation. In Section \ref{sec:other}, 
we derive (\ref{eqn:cubic}) and (\ref{eqn:MinusOne}) from the MV formula
and the cut-and-join equation. In Appendix A, we recall some well known facts
about Bernoulli numbers.

\bigskip

{\bf Acknowledgments}.
{\em This work was done in the Center of Mathematical Sciences, Zhejiang University.
The first and the third authors thank the Center for the hospitality and
the wonderful research environment}.

\section{Preliminaries}\label{sec:munford}

\subsection{Mumford's relations}

Let $\bE$ be the Hodge bundle over $\Mbar_{g,n}$,
and let $\lambda_i = c_i(\bE)$.
Define
\begin{align*}
c_t(\bE)& = \sum_{i=0}^g t^i\lambda_i, &
\Lambda_g^{\vee}(t) & = \sum_{i = 0}^g (-1)^i\lambda_it^{g-i}.
\end{align*}
Then we have
\begin{eqnarray*}
&& c_{-t}(\bE) = c_t(\bE^{\vee}) = (-t)^g \Lambda_g^{\vee}(\frac{1}{t}).
\end{eqnarray*}
Mumford's relations are given by:
\begin{eqnarray}
&& c_t(\bE)c_{-t}(\bE) = 1.
\end{eqnarray}
Equivalently,
\begin{eqnarray}
&& \Lambda_g^{\vee}(t)\Lambda_g^{\vee}(-t) = (-1)^gt^{2g}.
\end{eqnarray}

\subsection{Some consequences}
From the following well-known relation between Newton polynomials
and elementary symmetric polynomials (cf. e.g. \cite{Mac}):
$$\sum_{k \geq 1} p_k t^{k-1} = \frac{E_n'(-t)}{E_n(-t)},
$$
we get:
\begin{eqnarray} \label{eqn:Derivative1}
\sum_{n \geq 1} n! t^{n-1} \ch_n(\bE)
= \frac{c'_{-t}(\bE)}{c_{-t}(\bE)}= c_t(\bE)c'_{-t}(\bE).
\end{eqnarray}
It can be rewritten as
\begin{eqnarray*}
\sum_{n \geq 1} n! t^{n-1} \ch_n(\bE)
= \sum_{i=1}^g i \lambda_i (-t)^{i-1} \sum_{j=0}^g \lambda_jt^{j}.
\end{eqnarray*}
Hence $\ch_k(\bE) = 0$ for $k \geq 2g$, and
\begin{eqnarray}
&& n! \ch_n(\bE) = \sum_{i+j=n} (-1)^{i-1} i \lambda_i \lambda_j.
\end{eqnarray}
It is not hard to see that
\begin{eqnarray}
\ch_{2m}(\bE) & = & 0, \\
(2g-1)! \ch_{2g-1}(\bE) & = &(-1)^{g-1} \lambda_{g-1}\lambda_g, \\
(2g-3)! \ch_{2g-3}(\bE) & = &(-1)^{g-1}
(3\lambda_{g-3}\lambda_g - \lambda_{g-1}\lambda_{g-2}).
\end{eqnarray}

We will need the following results:

\begin{lemma} \label{lm:Derivative}
\begin{eqnarray}
&& \Lambda_g^{\vee}(1)(\Lambda_g^{\vee})'(-1)
= (-1)^{g-1}g +  \sum_{k \geq 1} k! (-1)^{k-1} \ch_k(\bE), \\
&& \left. \frac{d}{d\tau}\right|_{\tau=0} \left(\Lambda_g^{\vee}(1)
\Lambda_g^{\vee}(\tau)\Lambda_g^{\vee}(-\tau-1)\right)
= -\lambda_{g-1} + g \lambda_g \label{eqn:Derivative2} \\
&&\hspace{1.7in}  - \lambda_g \sum_{k \geq 1} k!(-1)^{k-1}\ch_k(\bE). \nonumber
\end{eqnarray}
In particular the degree $3g-3$ part of the left-hand side of
(\ref{eqn:Derivative2}) is
$$(-1)^{g-1}\lambda_g\lambda_{g-1}\lambda_{g-2}.$$
\end{lemma}

\begin{proof}
\begin{eqnarray*}
&& \Lambda_g^{\vee}(1)(\Lambda_g^{\vee})'(-1)
= \Lambda_g^{\vee}(1)\sum_{j=0}^g (-1)^j(g-j)t^{g-j-1}\lambda_j|_{t=-1} \\
& = &  (-1)^{g-1}g \Lambda_g^{\vee}(1)\sum_{j=0}^g \lambda_j
+ \Lambda_g^{\vee}(1) \cdot (-1)^g \sum_{j=0}^g j \lambda_j \\
& = & (-1)^{g-1} g + (-1)^gc_{-t}(\bE)c_t'(\bE)|_{t=1} \\
& = & (-1)^{g-1} g + (-1)^g \sum_{k \geq 1} k! (-1)^{k-1} \ch_k(\bE).
\end{eqnarray*}

\begin{eqnarray*}
&& \left. \frac{d}{d\tau}\right|_{\tau=0} \left(\Lambda_g^{\vee}(1)
\Lambda_g^{\vee}(\tau)\Lambda_g^{\vee}(-\tau-1)\right) \\
& = & \Lambda_g^{\vee}(1)\left. \frac{d}{d\tau}\right|_{\tau=0}
\Lambda_g^{\vee}(\tau) \cdot \Lambda_g^{\vee}(-1)
+ \Lambda_g^{\vee}(1) \Lambda_g^{\vee}(0) \cdot
\left. \frac{d}{d\tau}\right|_{\tau=0}\Lambda_g^{\vee}(-\tau-1) \\
& = & -\lambda_{g-1}
- (-1)^g \lambda_g \Lambda_g^{\vee}(1)(\Lambda^{\vee}_g)'(-1) \\
& = & -\lambda_{g-1} + g \lambda_g
- \lambda_g \sum_{k \geq 1} k!(-1)^{k-1}\ch_k(\bE).
\end{eqnarray*}
\end{proof}

\section{Mari\~no-Vafa Formula and ELSV Formula}
\label{sec:MVELSV}

In this section, we will explain close relationships between the
Mari\~{n}o-Vafa formula and the ELSV formula.

In Section \ref{sec:MV}, we recall the precise statement of the MV formula.
In Section \ref{sec:MVpf}, we describe our proof of the MV formula
\cite{LLZ1, LLZ} by the cut-and-join equation. In Section \ref{sec:ELSV},
we recall the ELSV formula relating Hurwitz numbers and Hodge integrals.
In Section \ref{sec:ELSVpf}, we describe a proof of the ELSV formula
completely parallel to the proof of the MV formula described in 
Section \ref{sec:MVpf}. In Section \ref{sec:burnside}, we will see
that a particular limit of the geometric and combinatorial sides of the MV
formula are the ELSV and Burnside formulas of Hurwitz numbers,
respectively.  

Both the MV formula and the ELSV formula  
relate the geometry of moduli spaces of Riemann surfaces encoded in Hodge integrals
to combinatorics of the representations of symmetric groups,
and hence to the theories of affine Kac-Moody Lie algebras and symmetric functions.

\subsection{Mari\~no-Vafa formula}\label{sec:MV}
The Mari\~{n}o-Vafa formula reads
\begin{equation}\label{eqn:MV}
\cC(\lambda;\tau;p)=R(\lambda;\tau;p)
\end{equation}
where $\cC(\lambda;\tau;p)$ is a generating function of Hodge integrals,
and $R(\lambda;\tau;p)$ is a generating function of representations
of symmetric group.

We first give the precise definition of $\cC(\lambda;\tau;p)$, the geometric
side of the MV formula.
For every partition $\mu = (\mu_1 \geq \cdots\geq \mu_{l(\mu)} > 0)$,
define
\begin{eqnarray*}
\cC_{g, \mu}(\tau)& = & - \frac{\sqrt{-1}^{|\mu|+l(\mu)}}{|\Aut(\mu)|}
 [\tau(\tau+1)]^{l(\mu)-1}
\prod_{i=1}^{l(\mu)}\frac{ \prod_{a=1}^{\mu_i-1} (\mu_i \tau+a)}{(\mu_i-1)!} \\
&& \cdot \int_{\Mbar_{g, l(\mu)}}
\frac{\Lambda^{\vee}_g(1)\Lambda^{\vee}_g(-\tau-1)\Lambda_g^{\vee}(\tau)}
{\prod_{i=1}^{l(\mu)}(1- \mu_i \psi_i)}, \\
\cC_{\mu}(\lambda; \tau) & = & \sum_{g \geq 0} \lambda^{2g-2+l(\mu)}\cC_{g, \mu}(\tau).
\end{eqnarray*}
Note that
$$
\int_{\Mbar_{0, l(\mu)}}
\frac{\Lambda^\vee_0(1)\Lambda^\vee_0(-\tau-1)\Lambda_0^\vee(\tau) }
{\prod_{i=1}^{l(\mu)} (1 -\mu_i \psi_i) }
=\int_{\Mbar_{0, l(\mu)}}
\frac{1}{\prod_{i=1}^{l(\mu)}(1 - \mu_i\psi_i)}
= |\mu|^{l(\mu)-3}
$$
for $l(\mu)\geq 3$, and we use this expression to extend the definition
to the case $l(\mu)<3$.
Introduce formal variables $p=(p_1,p_2,\ldots,p_n,\ldots)$, and define
$$
p_\mu=p_{\mu_1}\cdots p_{\mu_{l(\mu)} }
$$
for a partition $\mu=(\mu_1\geq \cdots \geq \mu_{l(\mu)}>0 )$.
Define the following generating series
\begin{eqnarray*}
\cC(\lambda; \tau; p) & = & \sum_{|\mu| \geq 1} \cC_{\mu}(\lambda;\tau)p_{\mu}.
\end{eqnarray*}

We next give the precise definition of $R(\lambda;\tau;p)$, the  combinatorial
side of the MV formula. For a partition $\mu$,
denote by $\chi_{\mu}$ the character of the irreducible representation
of $S_{|\mu|}$ indexed by $\mu$,
and by $C(\mu)$ the conjugacy class of $S_{|\mu|}$ indexed by $\mu$.
Define:
\begin{equation} \label{eqn:V}
\begin{split}
V_{\mu}(\lambda) = & \prod_{1 \leq a < b \leq l(\mu)}
\frac{\sin \left[(\mu_a - \mu_b + b - a)\lambda/2\right]}{\sin \left[(b-a)\lambda/2\right]} \\
& \cdot\frac{1}{\prod_{i=1}^{l(\nu)}\prod_{v=1}^{\mu_i} 2 \sin \left[(v-i+l(\mu))\lambda/2\right]}.
\end{split} \end{equation}
This has an interpretation in terms of quantum dimension
\cite{Mar-Vaf}. Let $q=e^{\lambda}$. Then
\begin{equation}
V_\mu(\lambda)=\frac{1}{\prod_{x\in\mu}(2\sinh((h(x)\lambda/2))}=
\frac{1}{\prod_{x\in\mu} (q^{h(x)/2}-q^{-h(x)/2})}= \frac{\dim_q R_\mu}{|\mu|!}.
\end{equation}
Define
\begin{eqnarray*}
& & R(\lambda; \tau; p) \\
&=&\log\left(\sum_{\mu}\left(\sum_{|\nu|=|\mu|}\frac{\chi_{\nu}(C(\mu))}{z_{\mu}}
e^{\sqrt{-1}(\tau+\frac{1}{2})\kappa_{\nu}\lambda/2} V_{\nu}(\lambda)
\right)p_\mu\right)\\
&=&\sum_{n \geq 1} \frac{(-1)^{n-1}}{n}  \sum_{\mu}
\left(\sum_{\cup_{i=1}^n \mu^i = \mu}
\prod_{i=1}^n \sum_{|\nu^i|=|\mu^i|} \frac{\chi_{\nu^i}(C(\mu^i))}{z_{\mu^i}}
e^{\sqrt{-1}(\tau+\frac{1}{2})\kappa_{\nu^i}\lambda/2} V_{\nu^i}(\lambda)
\right)p_\mu
\end{eqnarray*}

\subsection{Proof of Mari\~{n}o-Vafa formula by cut-and-join equation}
\label{sec:MVpf}

We briefly describe the proof given in \cite{LLZ1, LLZ}.

The combinatorial side $R(\lambda;\tau;p)$ of the MV formula
satisfies the follow cut-and-join equation:
\begin{equation}\label{eqn:MVcj}
\frac{\partial \Gamma}{\partial \tau}
= \frac{\sqrt{-1}\lambda}{2} \sum_{i, j\geq 1}
 \left(ijp_{i+j}\frac{\partial^2\Gamma}{\partial p_i\partial p_j}
+ ijp_{i+j}\frac{\partial \Gamma}{\partial p_i}\frac{\partial \Gamma}{\partial p_j}
+ (i+j)p_ip_j\frac{\partial \Gamma}{\partial p_{i+j}}\right).
\end{equation}

By functorial localization calculations  on moduli spaces 
$\Mbar_{g,0}(\bP^1,\mu)$ of relative stable maps,
the geometric side $\cC(\lambda;\tau;p)$ of the MV formula
also satisfies the cut-and-join equation (\ref{eqn:MVcj}).

Any solution $\Gamma(\lambda;\tau;p)$ to the cut-and-join equation
(\ref{eqn:MVcj}) is uniquely determined by its initial value
$\Gamma(\lambda;0;p)$, so the result follows from
$$
\cC(\lambda;0;p)=-\sum_{d\geq
  1}\frac{\sqrt{-1}^{d+1}p_d}
{2d\sin(d\lambda/2)}=R(\lambda;0;p)
$$
where the first equality follows from (\ref{eqn:sin})
and the second equality can be verified by combinatorics.

\subsection{ELSV formula}\label{sec:ELSV}

Given a partition $\mu$ of length $l(\mu)$,
denote by $H_{g_, \mu}$
the Hurwitz numbers of almost simple
Hurwitz covers of $\bP^1$ of ramification type $\mu$
by  connected genus $g$ Riemann surfaces.
The ELSV formula \cite{ELSV, Gra-Vak} states:
$$
H_{g, \mu} =(2g-2+|\mu|+l(\mu))! I_{g,\mu}
$$
where
$$
I_{g,\mu}=\frac{1}{|\Aut (\mu)|}
\prod_{i=1}^{l(\mu)} \frac{\mu_i^{\mu_i}}{\mu_i!}
\int_{\Mbar_{g, l(\mu)}} \frac{\Lambda_g^{\vee}(1)}{\prod_{i=1}^{l(\mu)} (1 - \mu_i \psi_i)}.
$$

Define generating functions
\begin{eqnarray*}
\Phi_{\mu}(\lambda)&=& \sum_{g \geq 0} H_{g, \mu}
\frac{\lambda^{2g-2+|\mu|+l(\mu)}}{(2g-2 + |\mu| + l(\mu))!}, \\
\Phi(\lambda; p) &=& \sum_{|\mu|\geq 1} \Phi_{\mu}(\lambda) p_{\mu}, \\
\Psi_{\mu}(\lambda) &=&\sum_{g \geq 0}I_{g,\mu}
\lambda^{2g-2+|\mu|+l(\mu)},\\
\Psi(\lambda; p) &=& \sum_{|\mu|\geq 1} \Psi_{\mu}(\lambda)p_{\mu}.
\end{eqnarray*}
In terms of generating functions,  the ELSV formula reads
\begin{equation}\label{eqn:ELSV}
\Psi(\lambda;p)=\Phi(\lambda;p).
\end{equation}

\subsection{Proof of ELSV formula by cut-and-join equation}
\label{sec:ELSVpf}

It was shown in \cite{Gou-Jac, Gou-Jac-Vai} that
$\Phi(\lambda;p)$ satisfies the following cut-and-join equation:
\begin{equation}\label{eqn:ELSVcj}
\frac{\partial \Theta}{\partial \lambda}
= \frac{1}{2} \sum_{i, j\geq 1} \left(ijp_{i+j}\frac{\partial^2\Theta}{\partial p_i\partial p_j}
+ ijp_{i+j}\frac{\partial \Theta}{\partial p_i}\frac{\partial \Theta}{\partial p_j}
+ (i+j)p_ip_j\frac{\partial \Theta}{\partial p_{i+j}}\right).
\end{equation}
Later (\ref{eqn:ELSVcj}) was reproved by sum formula of symplectic Gromov-Witten
invariants \cite{Li-Zha-Zhe, Ion-Par}. 

The calculations in Section 7 and Appendix A of \cite{LLZ} shows that
\begin{equation}\label{eqn:HI}
\tilde{H}_{g,\mu}= (2g-2+|\mu|+l(\mu))! I_{g,\mu}
\end{equation}
\begin{equation}\label{eqn:HIone}
\begin{split}
& \tilde{H}_{g,\mu}=(2g-3+|\mu|+l(\mu))!\left(
\sum_{\nu\in J(\mu)} I_{g,\nu}+\sum_{\nu\in C(\mu)} I_2(\nu) I_{g-1,\nu}\right.\\
&\makebox[3cm]{ } +\left.\sum_{g_1+g_2=g}\sum_{\nu^1\cup \nu^2\in C(\mu)} 
I_3(\nu^1,\nu^2)I_{g_1,\nu_1} I_{g_2,\nu_2}\right)
\end{split}
\end{equation}
where 
$$
\tilde{H}_{g,\mu}=\int_{[\Mbar_{g,0}(\bP^1,\mu)]^{\mathrm{vir} } }\mathrm{Br}^* H^r
$$ 
is some relative Gromov-Witten invariant of $(\bP^1,\infty)$, and 
$C(\mu), J(\mu), I_1, I_2,I_3$ are defined as in 
\cite{Li-Zha-Zhe}. So we  have
\begin{eqnarray*}
 &&(2g-2+|\mu|+l(\mu)) I_{g,\mu}\\
 &=&
\sum_{\nu\in J(\mu)} I_{g,\nu}+\sum_{\nu\in C(\mu)} I_2(\nu) I_{g-1,\nu}
+\sum_{g_1+g_2=g}\sum_{\nu^1\cup \nu^2\in C(\mu)} I_3(\nu^1,\nu^2)I_{g_1,\nu_1} I_{g_2,\nu_2},
\end{eqnarray*}
which is equivalent to the statement that the generating function $\Psi(\lambda;p)$ of $I_{g,\mu}$ also satisfies the cut-and-join equation (\ref{eqn:ELSVcj}).

Any solution $\Theta(\lambda;p)$ to the cut-and-join equation
(\ref{eqn:ELSVcj}) is uniquely determined
by its initial value $\Theta(0;p)$, so it remains to show that $\Psi(0;p)=\Phi(0;p)$.
Note that $2g-2+|\mu|+l(\mu)=0$ if and only if $g=0$ and $\mu=(1)$, so 
$$
\Psi(0;p)=H_{0,(1)}p_1,\ \ 
\Phi(0;p)=I_{0,(1)}p_1.
$$
It is easy to see that $H_{0,(1)}=I_{0,(1)}=1$, so
$$
\Psi(0;p)=\Phi(0;p).
$$

One can see geometrically that the relative Gromov-Witten invariant $\tilde{H}_{g,\mu}$
is equal to the Hurwitz number $H_{g,\mu}$. This together with (\ref{eqn:HI}) gives 
a proof of the ELSV formula presented in \cite[Section 7]{LLZ} in the spirit of \cite{Gra-Vak}.
Note that $\tilde{H}_{g,\mu}=H_{g,\mu}$  is not used in the proof described above.

\subsection{ELSV formula as limit of Mari\~{n}o-Vafa formula}
\label{sec:burnside}

By the Burnside formula, one easily gets
the following expression (see e.g. \cite{Zho}):
\begin{eqnarray*}
\Phi(\lambda;p)
&=&\log\left(\sum_\mu\left(\sum_{|\nu|=|\mu|}
\frac{\chi_{\nu}(\mu)}{z_{\mu}}
e^{\kappa_{\nu}\lambda/2} \frac{\dim R_{\nu}}{|\nu|!}\right) p_{\mu}.
\right)\\
& = & \sum_{ n \geq 1} \frac{(-1)^{n-1}}{n} \sum_{\mu}
\sum_{\cup_{i=1}^n \mu_i = \mu} \prod_{i=1}^n \sum_{|\nu_i|=|\mu_i|}
\frac{\chi_{\nu_i}(\mu_i)}{z_{\mu_i}}
e^{\kappa_{\nu_i}\lambda/2} \frac{\dim R_{\nu_i}}{|\nu_i|!} p_{\mu}.
\end{eqnarray*}
The ELSV formula reads
$$
\Psi(\lambda;p)=\Phi(\lambda;p)
$$
where the left hand side is a generating function of Hodge integrals $I_{g,\mu}$, and
the right hand side is a generating function of representations of symmetric groups.
So the ELSV formula and the MV formula are of the same type. 

Actually, the ELSV formula can be obtained by taking a particular
limit of the MV formula $\cC(\lambda;\tau;p)=R(\lambda;\tau;p)$.
More precisely, it is straightforward to check that
\begin{eqnarray*}
&& \lim_{\tau\to 0}\cC(\lambda\tau;\frac{1}{\tau};
(\lambda\tau)p_1,(\lambda\tau)^2 p_2,\cdots)\\
&=&\sum_{|\mu|\neq 0}
\sum_{g=0}^\infty \sqrt{-1}^{2g-2+|\mu|+l(\mu)} I_{g,\mu}\lambda^{2g-2+|\mu|+l(\mu)}p_\mu \\
&=& \Psi(\sqrt{-1}\lambda;p)
\end{eqnarray*}
and
\begin{eqnarray*}
&& \lim_{\tau\to 0}R(\lambda\tau;\frac{1}{\tau};
(\lambda\tau)p_1,(\lambda\tau)^2 p_2,\cdots)\\
&=&\log\left(\sum_{\mu}\left(\sum_{|\nu|=|\mu|}\frac{\chi_{\nu}(C(\mu))}{z_{\mu}}
e^{\sqrt{-1}\kappa_{\nu}\lambda/2}\lim_{t\to 0}
 (t^{|\nu|}V_{\nu}(t))\right) p_\mu\right)\\
&=&\log\left(\sum_{\mu}\left(\sum_{|\nu|=|\mu|}\frac{\chi_{\nu}(C(\mu))}{z_{\mu}}
e^{\sqrt{-1}\kappa_{\nu}\lambda/2} \frac{1}{\prod_{x\in\nu}h(x)}\right) p_\mu \right)\\
&=&\Phi(\sqrt{-1}\lambda;p)
\end{eqnarray*}
where we have used
$$
\frac{1}{\prod_{x\in\nu}h(x)}=\frac{\dim R_\nu}{|\nu|!}.
$$

In this limit, the cut-and-join equation (\ref{eqn:MVcj})
of $\cC(\lambda;\tau;p)$ and $R(\lambda;\tau;p)$ reduces 
to the cut-and-join equation (\ref{eqn:ELSVcj}) of 
$\Psi(\lambda;p)$ and $\Phi(\lambda;p)$, respectively.

\section{Mari\~{n}o-Vafa Formula and $\lambda_g$-Integrals}
\label{sec:MVlambda}
The $\lambda_g$ conjecture \cite{Get-P,Fab-Pan1} states
\begin{equation}\label{eqn:lambdagI}
\int_{\Mbar_{g, n}} \psi_1^{k_1} \cdots \psi_n^{k_n}\lambda_g =
\begin{pmatrix} 2g+n-3 \\ k_1, \dots, k_n\end{pmatrix} b_g
\end{equation}
where
\begin{equation}\label{eqn:bgI}
b_g=\left\{\begin{array}{ll}
1,& g=0,\\
\int_{\Mbar_{g,1}}\psi_1^{2g-2}\lambda_g, & g>0. 
\end{array}\right.
\end{equation}
The values of $b_g$ are given by \cite{Fab-Pan1}
\begin{equation}\label{eqn:sinI}
\sum_{g\geq 0} b_g t^{2g}=\frac{t/2}{\sin(t/2)}
=1+\sum_{g\geq 1} \frac{2^{2g-1}-1}{2^{2g-1}}\frac{|B_{2g}|}{(2g)!}t^{2g}
\end{equation}
where $B_{2g}$ are Bernoulli numbers. The $\lambda_g$-conjecture (\ref{eqn:lambdagI})
and the formula (\ref{eqn:bgI}) of $b_g$  give the following explicit formula
of $\lambda_g$-integrals: 
\begin{eqnarray} \label{eqn:lambdaB}
\int_{\Mbar_{g, n}} \psi_1^{k_1} \cdots \psi_n^{k_n}\lambda_g =
\begin{pmatrix} 2g+n-3 \\ k_1, \dots, k_n\end{pmatrix}
\frac{2^{2g-1}- 1}{2^{2g-1}} \frac{|B_{2g}|}{(2g)!}.
\end{eqnarray}

In this section, we will show that the formula (\ref{eqn:lambdaB})
of $\lambda_g$-integrals can be extracted from the Mari\~{n}o-Vafa 
formula, using the cut-and-join equation. Our proof of 
MV formula  relies only on  the formula (\ref{eqn:sinI})
of $b_g$ ($\lambda_g$ integrals on  $\Mbar_{g,1}$), so we obtain a new and simple proof of the 
$\lambda_g$-conjecture.

\subsection{A reformulation of (\ref{eqn:lambdaB}) }

We begin with the following reformulation.

\begin{lemma}\label{thm:formulate}
The formula (\ref{eqn:lambdaB}) of $\lambda_g$-integrals is equivalent to
\begin{eqnarray} \label{eqn:lambdag2}
&& \sum_{g \geq 0} \lambda^{2g}  \int_{\Mbar_{g,n}}
\frac{\lambda_g}{\prod_{i=1}^n (1- \mu_i\psi_i)}
= d^{n-3} \frac{d\lambda/2}{\sin (d\lambda/2)},
\end{eqnarray}
for all partitions of $d$.
\end{lemma}

\begin{proof}The left-hand side of (\ref{eqn:lambdag2}) is
\begin{eqnarray*}
&& \sum_{g \geq 0} \lambda^{2g} \int_{\Mbar_{g,n}}
\frac{\lambda_g}{\prod_{i=1}^n (1- \mu_i\psi_i)} \\
& = & \sum_{g \geq 0} \lambda^{2g} \sum_{k_1+\cdots + k_n = 2g-3+n}
\prod_{i=1}^n \mu_i^{k_i} \cdot
\int_{\Mbar_{g, n}} \lambda_g \prod_{i=1}^n \psi_i^{k_i}. \
\end{eqnarray*}
By (\ref{eqn:tsint}) in Appendix A
the right-hand side is:
\begin{eqnarray*}
&& d^{n-3} \left(1 + \sum_{g \geq 1} \frac{2^{2g-1}- 1}{2^{2g-1}} \frac{|B_{2g}|}{(2g)!}
(d\lambda)^{2g}\right) \\
& = & (\sum_i \mu_i)^{n-3}
+\sum_{g \geq 1} \lambda^{2g} \sum_{\sum_i k_i = 2g-3+n}
\prod_{i=1}^n \mu_i^{k_i} \cdot
\begin{pmatrix} 2g+n-3 \\ k_1, \dots, k_n\end{pmatrix}
\frac{2^{2g-1}- 1}{2^{2g-1}} \frac{|B_{2g}|}{(2g)!}.
\end{eqnarray*}
The Lemma is proved by comparing the coefficients.
\end{proof}

\subsection{Extracting (\ref{eqn:lambdaB}) from the Mari\~no-Vafa formula}
The following result has been proved in \cite{Zho} by  the cut-and-join equation.
\begin{theorem}
Write
$R(\lambda;\tau;p) = \sum_{\mu} R_{\mu}(\lambda;\tau)p_{\mu}$.
Then one has
\begin{equation} \label{eqn:Limit}
\begin{split}
& \lim_{\tau \to 0}  \lambda^{2-l(\mu)} \frac{1}{(\tau(\tau+1))^{l(\mu)-1}}
\prod_{i=1}^{l(\mu)} \frac{(\mu_i-1)!}{\prod_{j=1}^{\mu_i-1}(j+\mu_i\tau)}
\frac{\prod_j m_j(\mu)!}{\sqrt{-1}^{|\mu|+l(\mu)}}
R_{\mu}(\lambda;\tau) \\
= &  d^{l(\mu)-3} \cdot \frac{d\lambda/2}{\sin (d\lambda/2)}.
\end{split} \end{equation}
\end{theorem}

By the MV formula, the left-hand side of
(\ref{eqn:Limit}) is equal to
\begin{eqnarray*}
&& \sum_{g \geq 0} \lambda^{2g} \int_{\Mbar_{g, l(\mu)}}
\frac{\Lambda_g^{\vee}(1)\Lambda_g^{\vee}(0)\Lambda_g^{\vee}(-1)}{\prod_{i=1}^{l(\mu)}
(1 - \mu_i\psi_i)}
= \sum_{g \geq 0} \lambda^{2g} \int_{\Mbar_{g, l(\mu)}}
\frac{\lambda_g}{\prod_{i=1}^{l(\mu)} (1 - \mu_i\psi_i)}.
\end{eqnarray*}
Therefore,
we have established (\ref{eqn:lambdag2}) hence proved
the formula (\ref{eqn:lambdaB}) of $\lambda_g$ integrals.

\section{Derivation of Some Other Hodge Integral Identities}
\label{sec:other}

In this section we show how to derive from the Mari\~no-Vafa
formula and the cut-and-join equation the following formulas
proved in \cite{Fab-Pan1} by different methods:
\begin{eqnarray}
&& \int_{\Mbar_{g}} \lambda_{g-2}\lambda_{g-1}\lambda_g =
\frac{1}{2(2g-2)!} \frac{|B_{2g-2}|}{2g-2} \frac{|B_{2g}|}{2g}  \label{eqn:Cubic}\\
&& \int_{\Mbar_{g, 1}} \frac{\lambda_{g-1}}{1-\psi_1}
= b_g \sum_{i=1}^{2g-1} \frac{1}{i}
- \frac{1}{2} \sum_{\substack{g_1+g_2=g\\g_1,g_2 > 0}}
\frac{(2g_1-1)!(2g_2-1)!}{(2g-1)!}b_{g_1}b_{g_2} \label{eqn:g-1}
\end{eqnarray}

\subsection{The derivative}
We begin with the following special case of the Mari\~no-Vafa formula:
\begin{equation} \label{eqn:d}
\sum_{g \geq 0} \lambda^{2g} \int_{\Mbar_{g, 1}}
\frac{\Lambda^{\vee}_g(1)\Lambda^{\vee}_g(-\tau-1)\Lambda_g^{\vee}(\tau)}
{1- d \psi_i}
=  - \frac{\lambda}{\sqrt{-1}^{d+1}}
\frac{(d-1)!}{ \prod_{a=1}^{d-1} (d \tau+a)}
R_{(d)}(\lambda;\tau).
\end{equation}

\begin{lemma}
\begin{eqnarray} \label{eqn:Dd}
&& \left. \frac{d}{d \tau}\right|_{\tau=0}R_{(d)}(\lambda; \tau)
= \sum_{i+j=d, i \neq j}
\frac{-\sqrt{-1}^{d+1}\lambda}{8\sin (i\lambda/2)\sin (j\lambda/2)}.
\end{eqnarray}
\end{lemma}

\begin{proof}
This is an easy consequence of the cut-and-join equation:
\begin{eqnarray*}
\left. \frac{d}{d \tau}\right|_{\tau=0}R_{(d)}(\lambda; \tau)
& = &   \frac{\sqrt{-1}\lambda}{2} \sum_{i+j=d}
\left(ij R_{(i,j)}(\lambda;0)
+ ijR_{(i)}(\lambda; 0) R_{(j)}(\lambda;0)\right)\\
& = & \frac{\sqrt{-1}\lambda}{2} \sum_{i+j=d} i j
\frac{-\sqrt{-1}^{i+1}}{2i\sin (i\lambda/2)}
\frac{-\sqrt{-1}^{j+1}}{2j\sin (j\lambda/2)} \\
& = & \sum_{i+j=d}
\frac{-\sqrt{-1}^{d+1}\lambda}{8\sin (i\lambda/2)\sin (j\lambda/2)}.
\end{eqnarray*}
\end{proof}

\begin{corollary}
We have
\begin{equation} \label{eqn:DDD}
\begin{split}
& \sum_{g \geq 0} \lambda^{2g} \int_{\Mbar_{g, 1}}
\frac{\left. \frac{d}{d\tau}\right|_{\tau=0} \left(\Lambda_g^{\vee}(1)
\Lambda_g^{\vee}(\tau)\Lambda_g^{\vee}(-\tau-1)\right)}
{1- d \psi_1} \\
= & - \sum_{a=1}^{d-1} \frac{1}{a} \cdot \frac{d\lambda/2}{d\sin (d\lambda/2)}
+\sum_{i+j=d}
\frac{\lambda^2}{8\sin (i\lambda/2)\sin (j\lambda/2)}.
\end{split} \end{equation}
\end{corollary}

\begin{proof}
Take derivative in $\tau$ and set $\tau = 0$ on both sides of equation (\ref{eqn:d}).
By (\ref{eqn:Dd})
we get
\begin{eqnarray*}
&& \sum_{g \geq 0} \lambda^{2g} \int_{\Mbar_{g, 1}}
\frac{\left. \frac{d}{d\tau}\right|_{\tau=0} \left(\Lambda_g^{\vee}(1)
\Lambda_g^{\vee}(\tau)\Lambda_g^{\vee}(-\tau-1)\right)}
{1- d \psi_1} \\
& = & - \frac{\lambda}{\sqrt{-1}^{d+1}}
\left. \frac{d}{d \tau}\right|_{\tau=0}\frac{(d-1)!}{ \prod_{a=1}^{d-1} (d \tau+a)}\cdot
R_{(d)}(\lambda;0)
-  \frac{\lambda}{\sqrt{-1}^{d+1}} \left. \frac{d}{d \tau}\right|_{\tau=0}
R_{(d)}(\lambda;\tau) \\
& = & - \sum_{a=1}^{d-1} \frac{1}{a} \cdot \frac{d\lambda/2}{d\sin (d\lambda/2)}
+\sum_{i+j=d}
\frac{\lambda^2}{8\sin (i\lambda/2)\sin (j\lambda/2)}.
\end{eqnarray*}
\end{proof}

Now the left-hand side of (\ref{eqn:Dd}) is a polynomial in $d$
hence so must be the right-hand side.
If we find explicit expressions for the right-hand side,
then by comparing the coefficients,
we get Hodge integral identities.
This is how we prove (\ref{eqn:Cubic}) and (\ref{eqn:g-1}).

\subsection{The right-hand side}
We have
\begin{eqnarray*}
&& - \sum_{a=1}^{d-1} \frac{1}{a} \cdot \frac{d\lambda/2}{d\sin (d\lambda/2)}
= - \sum_{a=1}^{d-1} \frac{1}{a}
\sum_{g \geq 0} b_g  d^{2g-1}\lambda^{2g}
\end{eqnarray*}
hence the coefficient of $\lambda^{2g}$ is
\begin{eqnarray}
&&- \sum_{a=1}^{d-1} \frac{1}{a} \cdot
b_g d^{2g-1}.
\end{eqnarray}
This cancels with a similar term from the second term on the right-hand side of
(\ref{eqn:DDD}).

We also have
\begin{eqnarray*}
&&\sum_{i+j=d}
\frac{\lambda^2}{8\sin (i\lambda/2)\sin (j\lambda/2)} \\
& = & \sum_{i+j=d} \frac{1}{2ij}
 \sum_{g_1 \geq 0} b_{g_1} (i\lambda)^{2g_1}
\cdot \sum_{g_2 \geq 0}b_{g_2} (j\lambda)^{2g_2} \\
& = & \frac{1}{2} \sum_{g \geq 0} \lambda^{2g}
\sum_{g_1+g_2=g} b_{g_1}b_{g_2}  \sum_{i+j=d} i^{2g_1-1}j^{2g_2-1}.
\end{eqnarray*}

By (\ref{eqn:Prog}) in the Appendix we have for $g_1, g_2 > 0$, and $g_1+g_2=g$,
\begin{eqnarray*}
&&F_{g_1, g_2}(d)= \sum_{i+j=d} i^{2g_1-1}j^{2g_2-1}
= \sum_{i=1}^{d-1} i^{2g_1-1}(d-i)^{2g_2-1} \\
& = & \sum_{k=0}^{2g_2-1} (-1)^{2g_2-1-k}  \begin{pmatrix} 2g_2-1\\k\end{pmatrix}
d^{k} \sum_{i=1}^{d-1}i^{2g_1+2g_2-2-k} \\
& = & \sum_{k=0}^{2g_2-1} (-1)^{2g_2-1-k}   \begin{pmatrix} 2g_2-1\\k\end{pmatrix}
d^{k}  \sum_{l=0}^{2g-2-k}
 \frac{\begin{pmatrix} 2g-1-k \\ l \end{pmatrix}}{2g-1-k}B_ld^{2g-1-k-l} \\
& = & \sum_{k=0}^{2g_2-1}  \sum_{l=0}^{2g-2-k}
 \frac{(-1)^{2g_2-1-k} }{2g-1-k}  \begin{pmatrix} 2g_2-1\\k\end{pmatrix}
\begin{pmatrix} 2g-1-k \\ l \end{pmatrix}B_ld^{2g-1-l}.
\end{eqnarray*}
It is easy to see that the coefficient of $d$ in $F_{g_1, g_2}(d)$ receives
contribution only from the term with $k = 0$ and $l = 2g-2$,
hence it is
$$-B_{2g-2}.$$
The coefficient of $d^{2g-1}$ in $F_{g_1, g_2}(d)$ receives contributions from terms with
$l = 0$,
hence it is given by:
$$\sum_{k=0}^{2g_2-1}
\frac{(-1)^{2g_2-1-k}  }{2g-1-k}  \begin{pmatrix} 2g_2-1\\k\end{pmatrix}.$$

We can deal with the case of $g_1=0$ or $g_2 = 0$
in the same fashion.
\begin{eqnarray*}
&& F_{0, g}(d)
= \sum_{i+j=d} i^{-1}j^{2g-1}
= \sum_{i=1}^{d-1} i^{-1}(d-i)^{2g-1} \\
& = & \sum_{k=0}^{2g-1} (-1)^{2g-1-k}   \begin{pmatrix} 2g-1\\k\end{pmatrix}
d^{k} \sum_{i=1}^{d-1}i^{2g-2-k} \\
& = & \sum_{k=0}^{2g-2} (-1)^{2g-1-k}   \begin{pmatrix} 2g-1\\k\end{pmatrix}
d^{k}  \sum_{l=0}^{2g-k}
 \frac{\begin{pmatrix} 2g-1-k \\ l \end{pmatrix}}{2g-1-k}B_ld^{2g-1-k-l}
+ d^{2g-1} \sum_{i=1}^{d-1} \frac{1}{i} \\
& = & \sum_{k=0}^{2g-2}  \sum_{l=0}^{2g-2-k}
   \frac{(-1)^{2g-1-k}}{2g-1-k}  \begin{pmatrix} 2g-1\\k\end{pmatrix}
\begin{pmatrix} 2g-1-k \\ l \end{pmatrix}B_ld^{2g-1-l}
+ d^{2g-1} \sum_{i=1}^{d-1} \frac{1}{i}.
\end{eqnarray*}
The coefficient of $d$ in $F_{0, g}(d)$ or $F_{g. 0}(d)$ is
$$-B_{2g-2}.$$
The coefficient of $d^{2g-1}$ in $F_{0, g}(d)$ or $F_{g, 0}(d)$ is given by:
$$\sum_{k=0}^{2g-2}
\frac{(-1)^{2g-1-k}}{2g-1-k}  \begin{pmatrix} 2g-1\\k\end{pmatrix}
= \sum_{k=1}^{2g-1}
\frac{(-1)^i}{i}  \begin{pmatrix} 2g-1\\i\end{pmatrix}.$$

\subsection{Proof of (\ref{eqn:Cubic})}

By Lemma \ref{lm:Derivative},
the coefficient of $d\lambda^{2g}$ of the left-hand side of (\ref{eqn:DDD}) is:
\begin{equation} \label{eqn:LHS}
 (-1)^{g-1} \int_{\Mbar_{g, 1}}\psi_1
\lambda_{g} \lambda_{g-1}\lambda_{g-2} \\
= (-1)^{g-1}(2g-2) \int_{\Mbar_{g}}
\lambda_{g} \lambda_{g-1}\lambda_{g-2}.
 \end{equation}
By the above discussions,
the coefficient of $d\lambda^{2g}$
on the right-hand side of (\ref{eqn:DDD}) is
\begin{eqnarray*}
\frac{- B_{2g-2}}{2}\sum_{g_1+g_2=g}b_{g_1}
b_{g_2}
= \frac{- B_{2g-2}}{2} \cdot \frac{|B_{2g}|}{2g}
\frac{1}{(2g - 2)!}
\end{eqnarray*}
Comparing with (\ref{eqn:LHS}) we get
\begin{eqnarray*}
&& \int_{\Mbar_g} \lambda_g\lambda_{g-1}\lambda_{g-2}
= \frac{(-1)^{g}B_{2g-2}}{2(2g-2)} \cdot \frac{|B_{2g}|}{2g}
\frac{1}{(2g - 2)!}
\end{eqnarray*}
This is exactly (\ref{eqn:Cubic}).

\subsection{Proof of (\ref{eqn:g-1})}

The coefficient of $d^{2g-1}\lambda^{2g}$ on the left-hand side of (\ref{eqn:DDD}) is
\begin{eqnarray*}
&&- \int_{\Mbar_{g, 1}} \psi_1^{2g-1}\lambda_{g-1}
= - \int_{\Mbar_{g, 1}} \frac{\lambda_{g-1}}{1-\psi_1}.
\end{eqnarray*}
By the above discussions,
it is equal to
\begin{eqnarray*}
&& b_g \sum_{i=1}^{2g-1}
\frac{(-1)^i}{i}  \begin{pmatrix} 2g-1\\i\end{pmatrix}
+ \frac{1}{2} \sum_{\substack{g_1+g_2 = g\\ g_1, g_2 > 0}} b_{g_1}b_{g_2}
\sum_{k=0}^{2g_2-1}
\frac{(-1)^{2g_2-1-k}  }{2g-1-k}  \begin{pmatrix} 2g_2-1\\k\end{pmatrix}.
\end{eqnarray*}
Hence  (\ref{eqn:g-1}) is proved by the following:

\begin{lemma}
\begin{eqnarray*}
&& \sum_{k=0}^{2g_2-1}
\frac{(-1)^{2g_2-1-k}  }{2g-1-k}  \begin{pmatrix} 2g_2-1\\k\end{pmatrix}
= \frac{(2g_1-1)!(2g_2-1)!}{(2g-1)!}, \\
&& \sum_{i=1}^{2g-1}
\frac{(-1)^i}{i}  \begin{pmatrix} 2g-1\\i\end{pmatrix}
= -\sum_{i=1}^{2g-1} \frac{1}{i}.
\end{eqnarray*}
\end{lemma}

\begin{proof}
Let
$$f(x) = \sum_{i=1}^{2g-1} \frac{(-1)^i}{i}  \begin{pmatrix} 2g-1\\i\end{pmatrix}x^i.$$
Then we have
\begin{eqnarray*}
f'(x) & = & \sum_{i=1}^{2g-1} (-1)^i \begin{pmatrix} 2g-1\\i\end{pmatrix}x^{i-1}
= \frac{(1 - x)^{2g-1} - 1}{x}
= - \sum_{i=0}^{2g-2} (1-x)^i
\end{eqnarray*}
Hence we have
\begin{eqnarray*}
\sum_{i=1}^{2g-1} \frac{(-1)^i}{i} \begin{pmatrix} 2g-1\\i\end{pmatrix}
= \int_0^1 f'(x)dx = - \int_0^1 \sum_{i=0}^{2g-2} (1-x)^idx
= -\sum_{i=1}^{2g-1} \frac{1}{i}.
\end{eqnarray*}
Similarly,
let
$$g(x) = \sum_{k=0}^{2g_2-1}
\frac{(-1)^{2g_2-1-k}  }{2g-1-k}  \begin{pmatrix} 2g_2-1\\k\end{pmatrix}x^{2g-1-k}.$$
Then we have
\begin{eqnarray*}
g'(x) & = & \sum_{k=0}^{2g_2-1}
(-1)^{2g_2-1-k}  \begin{pmatrix} 2g_2-1\\k\end{pmatrix}x^{2g-2-k}
= x^{2g_1-1}(1-x)^{2g_2-1}.
\end{eqnarray*}
Hence we have by integrations by parts:
\begin{eqnarray*}
&& \sum_{k=0}^{2g_2-1}
\frac{(-1)^{2g_2-1-k}  }{2g-1-k}  \begin{pmatrix} 2g_2-1\\k\end{pmatrix} \\
& = &\int_0^1 g'(x)dx =  \int_0^1 x^{2g_1-1}(1-x)^{2g_2-1} dx \\
& = & \frac{(2g_2-1)}{2g_1} \int_0^1 x^{2g_1} (1 - x)^{2g_2-2} dx \\
& = & \frac{(2g_2-1)(2g_2-2)}{2g_1(2g_1+1)} \int_0^1 x^{2g_1+1} (1 - x)^{2g_2-3} dx  \\
& = & \cdots = \frac{(2g_2-1)!}{2g_1 (2g_1+1) \cdots (2g_1+2g_2-1)}
= \frac{(2g_1-1)!(2g_2-1)!}{(2g-1)!}.
\end{eqnarray*}
\end{proof}

\begin{appendix}

\section{Bernoulli numbers}
In this Appendix we recall some well known facts about Bernoulli numbers.
These numbers are defined by the following series expansion:
\begin{eqnarray} \label{eqn:B}
\frac{t}{e^t-1} = \sum_{m=0}^{\infty} B_m \frac{t^m}{m!}.
\end{eqnarray}
The first few terms are given by
\begin{align*}
B_0 & = 1, & B_1 & = - \frac{1}{2}, & B_2 & = \frac{1}{6}, & B_3 = 0.
\end{align*}
For odd $m > 1$,
$B_m = 0$,
and for even $m$,
the sign of $B_{2n}$ is $(-1)^{n-1}$.

\begin{lemma}
For $m > 0$,
\begin{eqnarray}
&& \sum_{k=0}^m \begin{pmatrix} m + 1 \\ k \end{pmatrix} B_k = 0. \label{eqn:B1}
\end{eqnarray}
\end{lemma}

\begin{proof}
Multiply both sides of (\ref{eqn:B}) by $e^t$.
The left-hand side becomes
\begin{eqnarray*}
&& e^t \frac{t}{e^t -1} = t + \frac{t}{e^t - 1}
= t + \sum_{m = 0}^{\infty} B_m \frac{t^m}{m!};
\end{eqnarray*}
the right-hand side becomes
\begin{eqnarray*}
&& \sum_{m = 0}^{\infty} B_m \frac{t^m}{m!} \cdot \sum_{n = 0}^{\infty} \frac{t^n}{n!}
= \sum_{m = 0}^{\infty} \sum_{k=0}^m B_k \frac{t^m}{k!(m-k)!}.
\end{eqnarray*}
Hence for $m > 1$,
we have
\begin{eqnarray*}
&& \frac{B_m}{m!} = \sum_{k=0}^m \frac{B_k}{k!(m-k)!}.
\end{eqnarray*}
(\ref{eqn:B1}) follows easily.
\end{proof}

\begin{lemma}
\begin{eqnarray}
&& \frac{t/2}{\sinh (t/2)}
= \sum_{m=0}^{\infty} \frac{1 - 2^{m-1}}{2^{m-1}}\frac{B_m}{m!} t^m,
\label{eqn:sinh} \\
&& \frac{t}{2}\coth (t/2)
= \sum_{n =0}^{\infty} B_{2n} \frac{t^{2n}}{(2n)!}. \label{eqn:coth}
\end{eqnarray}
\end{lemma}

\begin{proof}
These can be proved by easy algebraic manipulations as follows.
\begin{eqnarray*}
&& \frac{t/2}{\sinh (t/2)} = \frac{t}{e^t- 1} e^{t/2}
= 2 \frac{t/2}{e^{t/2} - 1} -  \frac{t}{e^t -1} \\
& = & 2\sum_{m=0}^{\infty} B_m \frac{(t/2)^m}{m!}
- \sum_{m=0}^{\infty} B_m \frac{t^m}{m!}
= \sum_{m=0}^{\infty} \frac{1 - 2^{m-1}}{2^{m-1}}\frac{B_m}{m!} t^m.
\end{eqnarray*}
\begin{eqnarray*}
\frac{t}{2}\coth (t/2) %% = \frac{t}{2} \frac{e^{t/2} + e^{-t/2}}{e^{t/2} - e^{-t/2}}
= \frac{t}{2} \frac{e^t+ 1}{e^t - 1}
= \frac{1}{2} + \frac{t}{e^t - 1}
= \sum_{n =0}^{\infty} B_{2n} \frac{t^{2n}}{(2n)!}.
\end{eqnarray*}
\end{proof}

\begin{corollary}
\begin{eqnarray}
&& \sum_{i+j=n} \frac{1-2^{1-2i}}{(2i)!}B_{2i} \cdot \frac{1 - 2^{1-2j}}{(2j)!}B_{2j}
= \frac{(1-2n)B_{2n}}{(2n)!}. \label{eqn:B2}
\end{eqnarray}
\end{corollary}

\begin{proof}
Apply the operator $t\frac{d}{dt}$ on both sides of (\ref{eqn:coth}):
\begin{eqnarray*}
&& \frac{t}{2} \coth\frac{t}{2} - \left(\frac{t/2}{\sinh (t/2)}\right)^2
=  \sum_{n=1}^{\infty} B_{2n}\frac{t^{2n}}{(2n-1)!}.
\end{eqnarray*}
Hence
\begin{eqnarray*}
&& \left(\frac{t/2}{\sinh (t/2)}\right)^2
=  \frac{t}{2} \coth\frac{t}{2} - \sum_{n=0}^{\infty} B_{2n}\frac{t^{2n}}{(2n-1)!}
= \sum_{n=0}^{\infty} (1-2n)B_{2n}\frac{t^{2n}}{(2n)!}.
\end{eqnarray*}
From this (\ref{eqn:B2}) easily follows.
\end{proof}

By changing $t$ to $\sqrt{-1}t$,
one gets from (\ref{eqn:sinh}) by recalling
$B_{2n} = (-1)^{n-1}|B_{2n|}$:
\begin{eqnarray} \label{eqn:tsint}
\frac{t/2}{\sin(t/2)}
= 1 + \sum_{g \geq 1} \frac{2^{2g-1}- 1}{2^{2g-1}} \frac{|B_{2g}|}{(2g)!} t^{2g}.
\end{eqnarray}
And (\ref{eqn:B2}) becomes
\begin{eqnarray} \label{eqn:BgBg}
\sum_{g_1+g_2=g}
\frac{2^{2g_1-1}- 1}{2^{2g_1-1}} \frac{|B_{2g_1}|}{(2g_1)!}
\frac{2^{2g_2-1}- 1}{2^{2g_2-1}} \frac{|B_{2g_2}|}{(2g_2)!} =
\frac{|B_{2g}|}{2g}
\frac{1}{(2g - 2)!}
\end{eqnarray}

Finally recall for any positive integer $m$,
\begin{eqnarray} \label{eqn:Prog}
&& \sum_{i=1}^{d-1} i^m = \sum_{k=0}^m
\frac{\begin{pmatrix} m +1 \\ k\end{pmatrix}}{m+1}B_kd^{m+1-k}.
\end{eqnarray}
\end{appendix}

\end{document}